\newfont{\bbb}{msbm10 scaled \magstep1}

 \newcommand{\proofhead}[1]{\par\pagebreak[1]\noindent{\bf#1.\ }}
 \newcommand{\pf}{\proofhead{Proof}}

\newcommand{\qed}{{\unskip\nolinebreak[1]\hspace{1.5em}\mbox{}\nolinebreak
    \hfill$\Box$\parfillskip=0pt\finalhyphendemerits=0\par\pagebreak[1]}}

\newtheorem{theorem}{Theorem}
\newtheorem{lemma}{Lemma}

\newtheorem{proposition}{Proposition}
\newtheorem{example}{Example}
\newtheorem{definition}{Definition}

\newtheorem{thm}{Theorem}
\newtheorem{prop}[thm]{Proposition}
\newtheorem{lem}{Lemma}
\newtheorem{rem}{Remark}
\newtheorem{cor}{Corollary}
\newtheorem{defn}{Definition}

\newcommand{\bed}{\begin{defn}}
\newcommand{\bet}{\begin{thm}}
\newcommand{\bep}{\begin{prop}}
\newcommand{\bel}{\begin{lem}}
\newcommand{\brk}{\begin{rem}}
\newcommand{\bdes}{\begin{description}}
\newcommand{\becor}{\begin{cor}}

\newcommand{\eed}{\end{defn}}
\newcommand{\eet}{\end{thm}}
\newcommand{\eep}{\end{prop}}
\newcommand{\eel}{\end{lem}}
\newcommand{\erk}{\end{rem}}
\newcommand{\edes}{\end{description}}

\newcommand{\eecor}{\end{cor}}

\newcommand{\be}{\begin{equation}}
\newcommand{\ee}{\end{equation}}
\newcommand{\beqa}{\begin{eqnarray*}}
\newcommand{\eeqa}{\end{eqnarray*}}
\newcommand{\beqn}{\begin{eqnarray}}
\newcommand{\eeqn}{\end{eqnarray}}
\newcommand{\eeq}{\end{equation}}
\newcommand{\beq}{\begin{equation}}

\newcommand{\non}{\noindent}

\documentstyle[11pt]{article}

\author{Horia C. Pop\\
University of Iowa}
\title{Row-Reducing the Quantum Determinant and Dieudonn\'e Determinant
}

\begin{document}
\date{}
\maketitle
\begin{abstract}

We  prove that row reducing a quantum matrix yields another  quantum matrix
for the same parameter $q$. This means that the elements of the new matrix 
satisfy the same 
relations as those of the original quantum matrix ring $M_q(n)$.
As a corollary, we can prove that the image of the 
quantum determinant in the abelianization of the total ring of quotients of 
$M_q(n)$ is equal to the Dieudonn\'e 
determinant of the quantum matrix. A similar result is proved for the
multiparameter quantum determinant. 
\end{abstract}

\section{Introduction}
\setcounter{thm}{0}

The representation theory of the quantum general linear group
is similar to that of  the commutative  
coordinate ring  of the general linear group. 

In this paper we show that row-reducing, an important method  in
linear algebra, does also work for 
quantum matrices and the quantum determinant.
 More precisely, in Theorem \ref{thm-row-red-qmatr}, we show that the 
matrix obtained from the original quantum matrix $Z$ by row reducing 
to clear all the elements, but one, in the first column, 
yields  a matrix that is again a quantum matrix, such that,  
the elements of the new matrix satisfy the same 
relations as those of the original quantum matrix ring $M_q(n)$. 
In Theorem \ref{thm-q-ind-step} we prove that 
the quantum determinant of the new (row-reduced) matrix, is equal to the 
quantum determinant of the original quantum matrix.
   
We work over an algebraically closed field $k$. The ring of quantum 
matrices $M_q(n)$ is an iterated Ore extension 
of $k$, hence it is an integral domain and it has a total ring of left 
(and right) quotients $D$.
We show in Corollary to Theorem \ref{thm-q-ind-step} that the image of the 
quantum determinant of the quantum matrix $Z$ in 
$D^{*}{}^{ab}$ is equal to the Dieudonn\'e determinant of the matrix $Z$.
This can also be established by using the remark 
following 3.1 and 3.3 in the paper [GR1] by Ghelfand and Retakh. 
The proofs of the above mentioned formulas are not provided, therefore,  
we shall 
give  a detailed account of our proof. On the other hand our  approach is 
quite different. It is  based on row-reducing a
 quantum matrix and on the fact that  the quantum determinant 
is invariant under row reducing over  $D$. Our proof in contrast 
to the approach in [GR1] and [GR2] has to do with a 
Grassmannian algebra. 

Although both are called determinants, the Dieudonn\'e determinant and the 
quantum determinant each represent two quite different approaches for the 
notion of a  noncommutative determinant. 
Basically, the Dieudonn\'e determinant comes 
from row reducing, while the quantum determinant comes from a Grassmannian 
algebra.  The major difficulty is caused by the fact that the quantum 
determinant is not, in general, multiplicative. 
 
The Dieudonn\'e determinant was defined in [Di] for an arbitrary square 
matrix with coefficients in a skew field; its values are cosets in 
$D^{*}{}^{ab}$. While one can provide a well-defined 
element in $D^{*}{}^{ab}$, we do not get a nice formula for it. To calculate 
the Dieudonn\'e determinant we use a procedure, similar to the one of 
producing the row-reduced echelon form of a matrix. 
 
The quantum determinant 
cannot be defined for an arbitrary matrix but only for 
a quantum matrix. For a quantum matrix we may consider  a 
quantum Grassmannian plane, which is a Frobenius quadratic algebra 
and has an associated quantum determinant, see [Ma], Ch.8. 
The value of the quantum determinant is given by an elegant polynomial formula.
It also has a nice formula for row or column expansion. In fact, the 
construction of the quantum determinant is possible because one has a 
categorical braiding which can be used to define the quantum Grassmannian 
algebra. Therefore, the quantum determinant corresponds to a setting which 
entails information not ordinarily  present when one just looks for a 
noncommutative determinant of a square matrix over an arbitrary skew field. 
This explains why the quantum determinant has properties that do not persist 
for the more general Dieudonn\'e determinant. 
 
For $n=2$ the expression of the quantum determinant is 
$det_q(Z)=ad-q^{-1}bc=ad-aca^{-1}b $,  therefore its projection 
in $D^{*ab}{}$, is equal to the Dieudonn\'e determinant. 
A direct calculation of the Dieudonn\'e determinant, for 
large $n$, does not look too nice anymore. 
%
We establish eventually the similar fact for an 
arbitrary $n$  by using an inductive argument.
 
The first two sections of this paper recall basic well known facts about the 
Dieudonn\'e determinant and the quantum determinant. Then we establish the 
Theorems \ref{thm-row-red-qmatr} and \ref{thm-q-ind-step} and as a corollary 
the fact that that 
the image of $det_q(Z)$ in $D^{*ab}$ equals the Dieudonn\'e determinant for 
the quantum matrix for any $n$. Finally we prove that a similar result 
holds for the multiparameter quantum determinant as defined in the paper by 
Artin, Schelter and  Tate [AST]. 
 
\medskip\  
 
\section{\bf The Dieudonn\'e determinant.} 
 
\smallskip\  
 
Let $D$ be a skew field and denote by $\pi $ the canonical surjection $%
D\longrightarrow D^{*}{}^{ab}=D^{*}/[D^{*},D^{*}].$ Recall some notations 
and definitions from [Dx]. Let $M_n(D)$ be the matrix units in the ring of 
square matrices of dimension $n$ over a skew field $D$, and denote by $%
e_{i,j}$ the matrix units in $M_n(D).$ 
 
As usual  we denote the transvections by $\tau _{i,j}(t)=I+te_{i,j},\ $ and we 
let $\ d_i(u)=I+(u-1)e_{i,i}.$ For an arbitrary permutation $\sigma \in S_n$ 
we associate the permutation matrix $P(\sigma )=(\delta _{i,\sigma j})\ $%
where $\delta $ denotes the Kronecker symbol. 
 
Multiplication on the left (right) by these matrices will perform the usual 
elementary row (column) transformations. 
 
Recall now that a matrix $A\in M_n(D)$ is right invertible if and only if it 
is left invertible, and this happens if and only if the rows (columns) of $%
A\ $are left (right) linearly independent. We denote the set of invertible 
matrices in $M_n(D)$ by $GL_n(D).$ 
 
The original design for the Dieudonn\'e determinant was to be a multiplicative 
mapping $\det _D:M_n(D)\longrightarrow D^{*}{}^{ab}$, equal to $1$ $\ $for 
any transvection $\tau _{i,j}(t)$. The price to pay for such a general 
definition is that its values are not in $D$ but  are cosets in $%
D^{*}{}^{ab}=D^{*}/[D^{*},D^{*}].\quad $It is difficult to handle this 
determinant because its values are not polynomials but rational functions in 
the coefficients and it hasn't any nice formula for row or column expansion. 
There is a substitute for the row expansion given by Theorem 4.5 , Ch.IV in 
[Ar], but it is too weak to be useful. 
 
\noindent We need the following result, see [Dx], Theorems 1 and 2 in Ch.19. 
 
\begin{theorem} [Bruhat normal form]
\label{Bruhat-normal-form}
 Let $D\;$be a skew field and $A\in M_n(D)\ $a matrix 
having a right inverse. Then there exists a decomposition$\ \quad $ 
 
\be 
A=TUP(\sigma )V  
\ee
 
$\ $where$\ T=\left(  
\begin{array}{ccc} 
1 & * & * \\   
& \ddots & * \\  
0 &  & 1  
\end{array} 
\right) ,\ U=diag(u_{1,}\ldots ,u_n)\ ,\ V=\left(  
\begin{array}{ccc} 
1 &  & 0 \\  
\,{}* & \ddots &  \\  
\,* & * & 1  
\end{array} 
\right) $ 
 
\smallskip\  
 
$\sigma $ a permutation, $P(\sigma )$ the permutation matrix corresponding to
 $\sigma$ and $U\ $and $\ \sigma \ $are unique with these properties. 
\end{theorem} 
 
\smallskip  
 
Now we may introduce the following: 
 
\begin{definition} 
\label{def-semi-det}
Let $A\in M_n(D),$and let $A=TUP(\sigma ), U=diag(u_{1,}\ldots ,u_n)\ ,$ with  
$u_{1,\cdots ,}u_n$ and $\sigma $ uniquely defined by the Bruhat normal form 
of $A.$ Define $\delta \epsilon \tau :M_n(D)\longrightarrow D,$ by: 
 
\qquad $\delta \epsilon \tau (A)=0$ if $A$ is not invertible and 
 
\qquad $\delta \epsilon \tau (A)=sign(\sigma )\ u_{1\cdot \;\ldots \cdot 
}\;u_n\ $if $A\in GL_n(D)$. 
\end{definition} 
 
\begin{definition} 
\label{def-ddet}
Let $\pi :D^{*}\longrightarrow D^{*}{}^{ab}=D^{*}/[D^{*},D^{*}]\ $ be the 
canonical projection. The Dieudonn\'e determinant of $A\ $is:  
$$ 
det_D(A)=\pi \ [\delta \epsilon \tau (A)].  
$$ 

\end{definition} 
 
Occasionally we denote the Dieudonn\'e determinant by $\left| A\right| _D.\ $ 
 
A direct consequence of this definition is that the Dieudonn\'e determinant 
has a number of remarkable properties which are similar to properties of a 
determinant in the commutative case: 
 
\begin{proposition}
\label{prop-ddet} 
1) $det_D(\tau _{i,j}(t))=1$ for any transvection $\tau $$%
_{i,j}(t)=I+te_{i,j}$ 
 
2)$\ det_D(P(\sigma ))=sign(\sigma )\ $for any permutation matrix $\sigma \ 
. $ 
 
3) The Dieudonn\'e determinant is multiplicative: $ det_D\left( A\cdot 
B\right) = det_D(A)\cdot det_D(B)$ 
\end{proposition} 
 
For the proofs see the original paper [Di], [Dx] Ch.19, 20 or [Ar] Ch.IV. 
 
The following two examples may help explain some of the reasons for 
taking cosets modulo commutators in D in defining the Dieudonn\'e 
determinant: 
 
\begin{example} 
\label{ex-prod-diag}
Let $U=\left(  
\begin{array}{ccc} 
u_1 &  & 0 \\   
& \ddots &  \\  
0 &  & u_n  
\end{array} 
\right) \quad $and $\quad V=\left(  
\begin{array}{ccc} 
v_1 &  & 0 \\   
& \ddots &  \\  
0 &  & v_n  
\end{array} 
\right) $ 
 
\smallskip\  
 
Then$\ \delta \epsilon \tau (U)=\ u_{1\cdot \;\ldots \cdot }\;u_n\ ,\delta 
\epsilon \tau (V)=v_{1\cdot \;\ldots \cdot }\;v_n\ $and $\delta \epsilon 
\tau (UV)=u_1v_{1\cdot \;\ldots \cdot }\;u_nv_n$ . 
 
This shows that $\delta \epsilon \tau $ is not multiplicative while $\pi 
(\delta \epsilon \tau (UV))=\pi (\delta \epsilon \tau (U))\cdot \pi (\delta 
\epsilon \tau (V)).$ 
 
Hence multiplicativity  will not hold unless we use the projection by $\pi $. 
\end{example} 
 
\begin{example} 
\label{ex-ord-diag}
Let $U=\left(  
\begin{array}{cc} 
u & 0 \\  
0 & v  
\end{array} 
\right) \ $and $V=\left(  
\begin{array}{cc} 
v & 0 \\  
0 & u  
\end{array} 
\right) $. Then V can be obtained from U by elementary row transformations, 
so that one can write $V=TU$ with$\ T$ a product of transvections. Here $%
\delta \epsilon \tau (U)=\ uv$ and $\delta \epsilon \tau (V)=vu$; they are 
not necessarily equal while $\pi (\delta \epsilon \tau (U))=\pi (\delta 
\epsilon \tau (V))$ 
\end{example} 
 
\noindent We calculate now the Dieudonn\'e determinant for the matrix $\left(  
\begin{array}{cc} 
a & b \\  
c & d  
\end{array} 
\right) $ when $a\neq 0.$ 
 
 
\begin{example}[The Dieudonn\'e determinant for a 2 by 2 matrix]
\label{ex-2dim-ddet}
\beqa\left(  
\begin{array}{cc} 
a & b \\  
c & d  
\end{array} 
\right) &=& \left(  
\begin{array}{cc} 
1 & 0 \\  
ca^{-1} & 1 
\end{array} 
\right) \left(  
\begin{array}{cc} 
a & b \\  
0 & d-ca^{-1}b  
\end{array} 
                \right) =\\ 
 \qquad &=& \left(  
\begin{array}{cc} 
1 & 0 \\  
ca^{-1} & 1 
\end{array} 
\right) \left(  
\begin{array}{cc} 
a & 0 \\  
0 & d-ca^{-1}b  
\end{array} 
\right) \left(  
\begin{array}{cc} 
1 & 0 \\  
0 & 1  
\end{array} 
\right) \left(  
\begin{array}{cc} 
1 & a^{-1}b \\  
0 & 1  
\end{array} 
\right) \eeqa 
 therefore 
 \beqn 
\left|  
\begin{array}{cc} 
a & b \\  
c & d  
\end{array} 
\right| _D=\pi (ad-aca^{-1}b). 
\eeqn 
 
Here we denoted the Dieudonn\'e determinant by $\left| A\right| _D.$ 
\end{example} 
 
\brk
\label{rem-non-prop-ddet}
Not all the properties for the commutative determinant will hold. For 
instance 
 
\begin{itemize} 
\item  there is no row expansion formula.
\item  the Dieudonn\'e determinant of the transpose of a matrix is not equal 
to the original determinant: 
\end{itemize} 
 \noindent For a noncommutative $D$ and two elements $a,b$ in $D$ such that $%
ab\neq ba,$the rows of the matrix $A=\left(  
\begin{array}{cc} 
1 & a \\  
b & ab  
\end{array} 
\right) $are left linearly independent hence the Dieudonn\'e determinant is 
non-zero. The transpose $A^t=\left(  
\begin{array}{cc} 
1 & b \\  
a & ab  
\end{array} 
\right) $has proportional rows. Hence $det_D(A)\neq det_D(A^t)=0$ see [Ar], 
Ch.IV, example following Thm.4.4. 
 
 \erk

\noindent We shall need the following result established in the original 
paper by Dieudonn\'e [Di], 1$^0$ in the proof of Thm.1. In fact Dieudonn\'e 
originally defined his noncommutative determinant by establishing first 
the following 
 
\begin{theorem} [Dieudonn\'e]
\label{thm-well-def-ddet}
Let $A=(a_{ij})$ and select a 
non-zero element in the first column, $a_{i1}\neq 0.\ $Let $A^{\prime 
}=(a_{ij}^{\prime })$ be the matrix obtained from $A$ by using elementary 
row operations to clear all elements in the first column except $a_{i1}$ 
i.e. subtract the $i$-th row multiplied on the left by $a_{j1}a_{i1}^{-1}$ 
from\ the j-th row. Denote by $\ A_{ij}$ the $(n-1)\times (n-1)$ matrix 
obtained by deleting the first column ant the $i$-th row from the original $%
A $. Then for any $i$ such that $a_{i1}\neq 0$ the product $\pi ($$%
a_{i1})\cdot \det _D(A_{i1})$ is the same, and it is equal to the Dieudonn\'e 
determinant $\det _D(A)$ of $A,$ i.e. $\det (A)=\pi (a_{i1})\cdot \det 
_D(A_{i1}).$ 
\end{theorem} 
 
\medskip\ \  
 
\section{\bf The Quantum determinant} 
 
\smallskip\  
 
The notion of quantum determinant was introduced in the theory of quantum 
groups, see [Dr], [Ma], [FRT] and [Gu]. The quantum determinant is a central 
element in the bialgebra of quantum matrices, and we shall see that it can 
be defined by a nice formula. It also has nice formulas for row or column 
expansion. Being a group-like element in the bialgebra of quantum matrices, 
it has a weak multiplicative property when applied to 
quantum matrices, such that their components commute pairwise. In that case 
the product of the two quantum matrices is again a point of $M_q(n).$ Keep 
also in mind that the quantum determinant can be defined for quantum 
matrices that are not square matrices (see Ch. 8, example 7, [Ma]). 
 
\begin{definition} 
\label{def-quadr-alg}
Let $A$ be a graded $k$-algebra $A=\oplus _{i\geq 0}A_i$ such that $A_0=k,$ 
dim$_k(A_1)<\infty ,$ $A$ is generated by $A_1$ as an algebra and the ideal 
of relations is generated by a set of homogeneous generators of degree $2.$ 
\end{definition} 
 Two very important examples of quadratic algebras are: 
 \begin{itemize} 
\item  The {\it quantum plane}: $A=k<x,y>/(yx-qxy)$ 
 \item  The {\it quantum Grassmannian:} $B=k<\xi ,\eta >/(\xi ^2,\eta ^2,\xi 
\eta -q\eta \xi )$. 
\end{itemize} 
 Another important exemple of the notion of quadratic algebra 
introduced by Manin is the quantum matrix bialgebra $M_q(2)$. Let $k$ be a 
field and $q\in k^{*}.\ $ 
 \begin{itemize} 
\item  For\ $n=2$\ the\ {\it quantum\ matrix\ bialgebra}\ is\ given\ by: 
 
$M_q(2)=k<a,b,c,d>$ subject to the following relations: 
 
$ab=q^{-1}ba,\ cd=q^{-1}dc$ 
 
$ac=q^{-1}ca,\ \ bd=q^{-1}db$ 
 
$bc=cb,\ ad-da=(q^{-1}-q)bc.$ 
\end{itemize} 
 The comultiple is defined as the dual of matrix multiplication namely: 
  
$\qquad \Delta (a)=a\otimes a+b\otimes c$ 
 
$\qquad \Delta (b)=a\otimes b+b\otimes d$ 
 
$\qquad \Delta (c)=c\otimes a+d\otimes c$ 
 
$\qquad \Delta (d)=c\otimes b+d\otimes d$ 
 
or in short: 
  
$$\qquad \Delta \left(  
\begin{array}{cc} 
a & b \\  
c & d  
\end{array} 
\right) =\left(  
\begin{array}{cc} 
a & b \\  
c & d  
\end{array} 
\right) \otimes \left(  
\begin{array}{cc} 
a & b \\  
c & d  
\end{array} 
\right) $$ 
 
\smallskip\  
 
The counit is determined by $\epsilon (a)=\epsilon (d)=1,\epsilon 
(b)=\epsilon (c)=0.$

Let A be an algebra over $k$ and let $\varphi :A\longrightarrow A$ be a
$k$-endomorphism. We call a linear endomorphism $\delta $ a $\varphi $%
-derivation of A if $\delta (ab)=\varphi (a)\delta (b)+a\delta (b),$ for all  
$a,b\in A$. For most of the following results one usually asks $\varphi $ to 
be injective, in our case it will be so granted that $\varphi $ is a 
$k$-automorphism. 
 
\begin{definition} 
\label{def-Ore-ext}
If $A\subseteq B,$then we say $B$ is a left Ore extension or skew polynomial 
extension of $A$ and write $B=A[X,\varphi ,\delta ]\ $if B is a ring of 
noncommutative polynomials in the variable$\ X$ with coefficients in $A\ $%
and the commutation relations between coefficients and variables are given 
by: $Xa=\varphi (a)X+\delta (a),$ for all $a\in A.$ 
\end{definition} 
 
We have then that $\ M_q(2)$ is an iterated Ore extension of $k$, as shown 
by the following tower of Ore extensions (see [Ka] Thm IV.4.10):  
\beqa k\hookrightarrow k[a,id,0]\hookrightarrow 
k[a][b,\varphi _{1,}0]\hookrightarrow k[a,b][c,\varphi _2,0]\hookrightarrow 
M_q(2)=k[a,b,c,d]=k[a,b,c][d,\varphi _3,\delta ]\eeqa 
 
\noindent  $\varphi _1:k[a]\ \rightarrow k[a]$ is a $k$-morphism such 
that $\varphi _1(a)=qa,$ 
 
\noindent $\varphi _2:k[\; a,b \;]\ \rightarrow \ k[\: a,b \: ]$ is a
 $k$-morphism such 
that $\varphi _2(a)=qa,\varphi _2(b)=b$ 
 
\noindent $\varphi _3:k[\: a,b,c\: ]\rightarrow k[\: a,b,c\: ]$ is a 
$k$-morphism with $\varphi 
_3(a)=a,\varphi _3(b)=qb,\varphi _3(c)=qc$ 
 
\noindent  $\delta \negthinspace :\negthinspace k[\: a,b,c \: ]\rightarrow 
k [\: a,b,c\: ] $ \ a $\varphi _3$-derivation defined by:
 $\delta (a)=(q-q^{-1})bc, \delta (b)=\delta (c)=0$.
 
When A is a domain and B is an Ore extension then B is also a domain. For 
this reason $M_q(2)$ is a domain. When A is left noetherian, then the Ore 
extension B is also left noetherian ([Ka] Proposition I.8.3), therefore the 
left Ore condition is satisfied (see [Co] proposition 1.3.2) for the 
multiplicative system of non-zero elements and one can construct a total 
field of fractions of $M_q(2),\ $that will be denoted in the sequel by $D.$ 
 
Let $E=k<a,b,c,d>$ with $a,b,c,d$ subject to the relations in $M_q(2)$ and
 $R=E\otimes _kB.$ It means $R$ is generated over $\ k$ 
by $a,b,c,d$ and by $\xi $ and $\eta $ 
(subject to the relations in $B),$ and any of $a,b,c,d$ commutes with 
any of $%
\xi $ and $\eta $ . Let us denote : 
 
$x^{\prime }=ax+by,\quad y^{\prime }=cx+dy$ \quad and\quad 
$\xi ^{\prime 
}=a\xi +b\eta ,\quad \eta ^{\prime }=c\xi +d\eta $ 
 
We say $\ x^{\prime },y^{\prime }$ is a point of the quantum plane if $%
y^{\prime }x^{\prime }=q\,x^{\prime }y^{\prime }$ and $\xi ^{\prime },\eta 
^{\prime }$ is a point of the quantum Grassmannian if $\xi ^{\prime 2}=0,\eta 
^{\prime 2}=0,$ and $\xi ^{\prime }\eta ^{\prime }=q\,\eta ^{\prime }\xi 
^{\prime }.$ 
 
Kobyzev noticed, see [Ma], Ch1, that $a,b,c$ and $d$ satisfy the relations 
defining $M_q(2),$ if and only if $x^{\prime },y^{\prime }$ is a point of 
the quantum plane and $\xi ^{\prime },\eta ^{\prime }$ is a point of the 
quantum Grassmannian. 
 
The quantum determinant for $n=2$ is defined to be:\ $det_q=ad-q^{-1}bc$ , 
this is just the unique element $det_q(Z)$ in $M_q(2)$ such that $\xi 
^{\prime }\eta ^{\prime }$ $=$ $det_q(Z)\cdot \xi \eta $ in $R.$ 
 
Using now the commutation relations in $M_q(2)$ we get $%
det_q(Z)=ad-aca^{-1}b .$ 
 
If we think of the quantum determinant as being an element in $D$ associated 
to the matrix $\left(  
\begin{array}{cc} 
a & b \\  
c & d  
\end{array} 
\right) $ then the image of $det_q(Z)$ in $D^{*}{}^{ab}$ is equal to the 
value of the Dieudonn\'e determinant for this matrix. 
 
On the other side direct calculation shows that 
the quantum determinant $det_q(Z)$ is a group-like element and 
$\ $commutes with $a,b,c,d$. Therefore  
$det_q$ is central in $M_q(2)$ . 
 
For an arbitrary $n$ define 
 
\begin{itemize} 
\item  the {\it quantum plane}:$\ A=k<x_{1,}\cdots ,x_n>$ subject to:$\quad 
x_jx_i-qx_ix_j=0$ for all $i<j\quad $
 
\item  the {\it quantum Grassmannian plane}: $B=k<\xi _{1,}\cdots ,\xi _n>$ 
subject to : $\xi _i^2=0$ for all $i$, and $\xi _i\xi _j-q^{-1}\xi _j\xi 
_i=0 $ for all $i<j$ 
 
\item  the bialgebra of {\it quantum matrices} to be the noncommutative ring 
\end{itemize} 
 
$M_q(n)=k<z_{11,}z_{12,}\cdots ,z_{nn}>$ \ subject to the following relations: 
 
 \beqa 
{\cal R}el_q  \left \{ \begin{array}{l}
 \mbox{\it Row relations:\ } z_{il}\ z_{ik}=q\ z_{ik}\ z_{il} 
\mbox{\it \  for all } k<l \label{mq-rel-row}\\
 \mbox{\it Column relations:\ } z_{jk}\ z_{ik}=q\ z_{ik}\ z_{jk}
\mbox{\it \ for all } i<j \label{mq-rel-col}\\
 \mbox{\it Secondary diagonal relations:\ } z_{il}\ z_{jk}=z_{jk}\ z_{il}
 \mbox{\it \ for all }  i<j \mbox{\it , } k<l   \label{mq-rel-sdiag}\\
 \mbox{\it Main diagonal relations:\ } z_{jl}\ z_{ik}-z_{ik}\: z_{jl}=
(q^{-1}-q)\: z_{il}\ z_{jk} \: \mbox{\it \ for all } i<j \mbox{\it ,} k<l. 
\label{ma-rel-mdiag}
\end{array}
\right.
\eeqa
 
\smallskip\  
 
\noindent We denote by $Z$ the matrix $(z_{ij})$ such that $z_{ij}$ satisfy 
the relations ${\cal R}el_q$ of $M_q(n).$ 
 
\noindent The comultiplication on $M_q(n)$ is dual to matrix multiplication, 
i.e., $\ \Delta (z_{ij})=\sum\limits_kz_{ik}\otimes z_{kj}$ and the counit is 
defined by $\epsilon (z_{ij})=\delta _{ij,}$ both extended multiplicatively 
to all $M_q(n).$ Then $\ M_q(n)$ is a bialgebra. 
 
Notice that one way of looking at these relations is asking whether for any 
selection of 2 rows and 2 columns the 4 elements at the intersection of the 
selected rows and columns satisfy the relations demanded before for $a,b,c,d$ 
in the case $n=2.$ 
 
\noindent Denote again $E=k<z_{ij}>,R=E\otimes _kB$ , the observation made 
by Kobyzev still holds: 
 
\begin{proposition} 
\label{Kobyzev-char-qmatr}
Let  
$$ 
x_i^{\prime }=\sum_jz_{ij}x_j\quad \mbox{and\quad }\xi _i^{\prime 
}=\sum_jz_{ij}\xi _j.  
$$ 
Then $z_{ij}$ satisfy the defining equations of $M_q(n)\ $if and only if $%
x_i^{\prime }$ (respectively $\xi _i^{\prime }$) satisfy the defining 
relations of the quantum plane (respectively the quantum Grassmannian 
plane$).$ 
\end{proposition} 
 
Now we follow Manin [Ma], Ch8 to define the quantum determinant

\begin{definition} 
\label{def-Frob-QA}
A quadratic algebra $A$ is called a Frobenius algebra of dimension $d$ if $%
dim(A_d)=1$, $A_i=0$ for $i>d$ and for all $j$ , $0\leq j\leq d$ the 
multiplication map $m:A_j\otimes A_{d-j}\longrightarrow A_d$ is a perfect 
duality. If in addition $dimA=\left(  
\begin{array}{c} 
d \\  
i  
\end{array} 
\right) $, call A a quantum Grassmannian algebra. 
\end{definition} 
 
\begin{definition} 
\label{def-man-Frob-qdet}
Let A be a Frobenius algebra and let E be a bialgebra coacting on A by $%
\delta :A\longrightarrow E\otimes A$ and such that $\delta (A_1)\subseteq 
E\otimes A_1.$ Then by induction, for any $j$ one has $\delta (A_j)\subseteq 
E\otimes A_{j,}$ and in particular since $dim(A_d)=1,$ there is an element 
called the quantum determinant of the coaction $D=DET(\delta )\in E$ such 
that for any a$\in A_d$ one has:  
$$ 
\delta (a)=DET(\delta )\otimes a  
$$ 
\end{definition}

\noindent An immediate consequence of coassociativity is that $DET(\delta )$ 
is a group-like element:  
$$ 
\Delta (DET(\delta ))=DET(\delta )\otimes DET(\delta ).  
$$ 
 
It is easy to see now, cf. [Ma] Ch8 Example 6 that: 
 
\begin{proposition} 
\label{qdet-formula}
The quantum Grassmannian plane is a Frobenius quadratic algebra, in fact even 
a quantum Grassmannian algebra. The quantum determinant in this case is
 given by 
the formula: 
  
\beq \qquad \left| Z\right| _q=det_q(Z)=\sum_{\sigma \in S_n} 
(-q^{-1})^{l(\sigma )}\ z_{1\sigma (1)\ }z_{2\sigma (2)\ }\cdots \ 
z_{n\sigma (n)} \eeq 
  
where $l(\sigma )$ is the length of $\sigma ,$ equal to the number of 
inversions in $\sigma $. 
\end{proposition} 
 
\begin{proposition} 
\label{grassm-def-qdet}
Let $E=k<z_{ij}>,R=E\otimes _kB,$ and $\xi _i^{\prime }=\sum_jz_{ij}\xi _j.$ 
 Then in $R$ we have: 
 
$$ 
\xi _1^{\prime }\xi _2^{\prime }\cdots \xi _n^{\prime }=det_q(Z)\ \xi _1\xi 
_2\cdots \xi _n  
$$ 
\end{proposition} 
 
An immediate consequence of the definition of the quantum determinant is 
the Laplace expansion formulas for rows and columns expansion. Denote by 
$Z_{ji}$ the quantum determinant of the $(n-1)$ by $(n-1)$ matrix obtained by 
removing the $j$-th row and the $i$-th column of $Z$ . 
$\tilde Z=((-q)^{j-i}Z_{ji})$ is called the $q$-cofactor matrix.
The following is proved in [FRT] Thm 4, see also [Tk] Proposition 2.3 and 
[PW]. 
 
\begin{proposition} 
\label{prop-qmat-inv}
If Z is a quantum matrix then $Z\tilde Z=\tilde ZZ=det_q(Z)I$ . Consequently  
$\quad det_q(Z)$ is central and the following row expansion formulas hold: 
\end{proposition} 
 
$$ 
det_q(Z)\ =\sum_j(-q)^{i-j}z_{ij}Z_{ij}=\sum_i(-q)^{i-j}z_{ij}Z_{ij}.  
$$ 
There is an algebra automorphism $\tau :M_q(n)\longrightarrow M_q(n),$ 
called the transposition and defined by extending multiplicatively from $%
\tau (z_{ij})=z_{ji}.$ It is known that the quantum determinant is invariant 
under $\tau .$ 
 
We remark now that exactly as in the case $n=2$ shown before the ring of
 quantum matrices $ M_q(n)=k<z_{11,}z_{12,}\cdots ,z_{nn}>$
 is an iterated Ore extension of $k$, 
therefore we can establish inductively that it is a left and right 
noetherian domain. Hence $M_q(n)$ is a domain and satisfies the Ore 
condition for the multiplicative system of non-zero elements. From now on we 
shall denote by $D$ the skew field that is the total ring of fractions of $%
M_q(n)=k<z_{11,}z_{12,}\cdots ,z_{nn}>.$ 
 
\begin{definition} 
\label{def-qmat-point} 
Let $R$ be a $k$-algebra. We call an $R$-point of $M_q(n)$ , (respectively 
of the quantum plane , or quantum Grassmannian plane) any $k$ -algebra 
morphism in $A\lg {}_k(M_q(n),R)$ (respectively in $A\lg {}_k(A,R)$ or $A\lg 
{}_k(B,R)$) 
\end{definition} 

Any $k$-algebra morphism $\Phi :$ $M_q(n)\longrightarrow R$ is uniquely 
determined by a matrix with entries in $R$ given by $r=(r_{ij})=\Phi 
(z_{ij}) $ such that the elements $r_{ij}$ satisfy the relations ${\cal R}%
el_q{\cal \ }$for quantum matrices. Given such a matrix we shall denote by $%
\Phi _r$ the $R$ point of $M_q(n)$ uniquely determined by $r,$ and we shall 
say that $r$ is a quantum matrix . 
 
\begin{definition}
\label{def-qdet-eval} 
The quantum determinant evaluated at an $R$ point $r$ of $M_q(n)$ is 
 
\smallskip\  
 
$\qquad \left| r\right| _q=\Phi _r(\det _q(X))=\sum_{\sigma \in S_n}$ $%
(-q^{-1})^{l(\sigma )}\ r_{1\sigma (1)\ }r_{2\sigma (2)\ }\cdots \ 
r_{n\sigma (n)}.$ 
 
\smallskip\  
\end{definition} 
 
\section{\bf The main results}
 
\smallskip\  
 
Let us now define a new matrix $Z^{\prime }=(z_{ij}^{\prime })$ by : 
 \beqn z_{1j}^{\prime }&=&z_{1j} \mbox{ for any } j,z_{i1}^{\prime }=0
 \mbox{ for } i=2,\ldots ,n \label{z'-1-st-row}\\
 z_{ij}^{\prime }&=&z_{ij}-z_{i1}z_{11}^{-1}z_{1j} \mbox{for } i=2,\ldots ,n
\mbox{ and} j=2,\ldots ,n. \label{z'-i-th-row}\eeqn 
 Note that $Z^{\prime }$ is obtained from $Z$ by clearing up all 
the positions in the first column except the first one by use of elementary 
row transformations. 
 
Let $Z^{\prime \prime }$ be the $(n-1)\times (n-1)$ matrix obtained by 
deleting the first row and the first column of $Z^{\prime }.$ 
 
\begin{theorem}[Row-reducing a quantum matrix]
\label{thm-row-red-qmatr}
The matrix $Z^{\prime }$ is a quantum matrix, its elements satisfy the 
defining relations ${\cal R}el_q{\cal \ }$of $M_q(n).$ 
\end{theorem} 
 
{\pf} 
For a selection of first row with any other row and any two columns 
different from the first, we may check directly that the relations are 
satisfied. An observation that makes calculations easier is that $%
z_{ij}^{\prime }=z_{11}^{-1}\left|  
\begin{array}{cc} 
z_{11} & z_{1i} \\  
z_{j1} & z_{ji}  
\end{array} 
\right| _q,$ now for $a=z_{11},b=z_{1i},c=z_{j1},d=z_{ji}$ one  may use 
the fact that the quantum determinant $ad-q^{-1}bc$ commutes with $a,b,c,d.$ 
 
For the other selections direct calculations will prove that if $%
A=z_{ik}^{\prime }\ ,B=z_{il}^{\prime }\ ,C=z_{jk}^{\prime }\ 
,D=z_{jl}^{\prime }$ $\ $ and $i<j,\ k<l$ then $A,B,C$ and D satisfy the 
relations for quantum matrices. 
of general $i,j,k,l$\ none of them $=1$) works in exactly the same way. 
All the proofs work  in a similar way: we express all monomials in 
terms of a Poincar\'e-Birkhoff-Witt basis of monomials in $z_{ij}^{\prime }s$ 
with indices in increasing lexicographic order by using the commutation 
relations and look for cancellations. 
 
\  
 
\begin{lemma} 
\label{lem-bc=cb}
$BC=CB$ 
\end{lemma} 
 
{\pf} 
We express $BC-CB$ in terms of a basis of monomials 
in $z_{ij}^{\prime }s$ with indices in increasing lexicographic order using 
the commutation relations 
 
$%
BC-CB=(z_{il}-z_{i1}z_{11}^{-1}z_{1l})(z_{jk}-z_{j1}z_{11}^{-1}z_{1k})-
(z_{jk}-z_{j1}z_{11}^{-1}z_{1k})(z_{il}-z_{i1}z_{11}^{-1}z_{1l})=  $ 
 
$\underbrace{z_{il}z_{jk}}_{\frame{1}} - \underbrace{z_{il} z_{j1} 
z_{11}^{-1} z_{1k}}_{\frame{2}} -
\underbrace{z_{i1}z_{11}^{-1}z_{1l}z_{jk}}_{\frame{3}} + 
\underbrace{ z_{i1}z_{11}^{-1} z_{1l}z_{j1}z_{11}^{-1} z_{1k}}_{\frame{4}} +$ 
 
 $ - \underbrace{z_{jk}z_{il}}_{\frame{1'}} +
 \underbrace{z_{j1}z_{11}^{-1}z_{1k}z_{il}}_{\frame{2'}} +
\underbrace{z_{jk}z_{i1}z_{11}^{-1}z_{1l}}_{\frame{3'}} -
\underbrace{z_{j1}z_{11}^{-1}z_{1k}z_{i1}z_{11}^{-1}z_{1l}}_{\frame{4'}}$ 
 
$\frame{2'}=q^{-1}z_{11}^{-1}z_{1k}z_{il}z_{j1}$ 
 
$\frame{3}=q^{-1}z_{11}^{-1}z_{1l}z_{i1}z_{jk}$ 
 
$\frame{2}=q^{-1}(z_{il}z_{11}^{-1})z_{1k}z_{j1}=($ use now the formula $%
da^{-1}-a^{-1}d=(q^{-1}-q)a^{-1}bca^{-1}=$ 
 
$=(q^{-1}-q)q^{-2}a^{-2}bc$ $)$ 
 
$%
=q^{-1}z_{11}^{-1}(z_{il}z_{1k})z_{j1}+
(q^{-1}-q)q^{-2}z_{11}^{-2}z_{1k}z_{1l}z_{i1}z_{j1}=(  
$use now $ad-da=(q^{-1}-q)bc$) 
 
$%
=q^{-1}z_{11}^{-1}z_{1k}z_{il}z_{j1}-(q^{-1}-q)q^{-1}z_{11}^{-1}z_{1l}
z_{ik}z_{j1}+(q^{-1}-q)q^{-2}z_{11}^{-2}z_{1k}z_{1l}z_{i1}z_{j1}  
$ 
 
$\frame{3'}%
=q^{-1}(z_{jk}z_{11}^{-1})z_{1l}z_{i1}=q^{-1}z_{11}^{-1}z_{1l}(z_{jk}z_{i1})+
(q^{-1}-q)q^{-2}z_{11}^{-2}z_{1k}z_{1l}z_{i1}z_{j1}=  
$ 
 
$%
=q^{-1}z_{11}^{-1}z_{1l}z_{i1}z_{jk}-(q^{-1}-q)q^{-1}z_{11}^{-1}z_{1l}
z_{ik}z_{j1}+(q^{-1}-q)q^{-2}z_{11}^{-2}z_{1k}z_{1l}z_{i1}z_{j1}  
$ \
 
\smallskip\  
 
Now one can see that all like terms cancel. \qed
 
\begin{lemma} 
\label{lem-ab=qab}
$BA=qAB$ 
\end{lemma} 
 
{\pf}
 
$BA-qAB=$ 
 
$=(z_{il}-z_{i1}z_{11}^{-1}z_{1l})(z_{ik}-z_{i1}z_{11}^{-1}z_{1k}) 
 -q (z_{ik}-z_{i1}z_{11}^{-1}z_{1k})(z_{il}-z_{i1}z_{11}^{-1}z_{1l})=  $ 
 
$\underbrace{z_{il}z_{ik}}_{\frame{1}} - \underbrace{z_{il}z_{i1}z_{11}^{-1}
z_{1k}}_{\frame{2}} - \underbrace{z_{i1}z_{11}^{-1}z_{1l}z_{ik}}_{\frame{3}} + 
\underbrace{z_{i1}z_{11}^{-1}z_{1l}z_{i1}z_{11}^{-1}z_{1k}}_{\frame{4}} +$ 
 
$- \underbrace{q z_{ik}z_{il}}_{\frame{1'}}  + 
q \underbrace{z_{i1}z_{11}^{-1}z_{1k}z_{il}}_{\frame{2'}} +
q \underbrace{z_{ik}z_{i1}z_{11}^{-1}z_{1l}}_{\frame{3'}} -
q \underbrace{z_{i1}z_{11}^{-1}z_{1k}z_{i1}z_{11}^{-1}z_{1l}}_{\frame{4'}}$ 
 
Now $\frame{1}=\frame{1'}$ and $\frame{4}=\frame{4'}$ and: 
 
$\frame{2'}=z_{11}^{-1}z_{1k}z_{i1}z_{il}$ 
 
$\frame{3}=q^{-1}z_{11}^{-1}z_{1l}z_{i1}z_{ik}$ 
 
$\frame{2}%
=q^{-1}(z_{il}z_{11}^{-1})z_{i1}z_{1k}=q^{-1}z_{11}^{-1}(z_{il}z_{1k})z_{i1}-
(q^{-1}-q)q^{-2}z_{11}^{-2}z_{1l}z_{i1}z_{1k}z_{i1}=  
$ 
 
$=q^{-1}z_{11}^{-1}z_{1k}z_{il}z_{i1}\ 
-(q^{-1}-q)z_{11}^{-1}z_{1l}z_{i1}z_{ik}-(q^{-1}-q)q^{-2}z_{11}^{-2}z_{1k}
z_{1l}z_{i1}^2  
$ 
 
$\frame{3'}%
=(z_{ik}z_{11}^{-1})z_{i1}z_{1l}=qz_{11}^{-1}z_{1l}z_{i1}z_{ik}+(q^{-1}-q)
q^{-2}z_{11}^{-2}z_{1k}z_{1l}z_{i1}^2  
$  
 
\smallskip\  
 
The cancellations are now clear. \qed
 
\begin{lemma} 
$CA=qAC$ 
\end{lemma} 
 \label{lem-ca=qac}
{\pf}\smallskip\  
 
$CA-qAC=$ 
 
$%
=(z_{jk}-z_{j1}z_{11}^{-1}z_{1k})(z_{ik}-z_{i1}z_{11}^{-1}z_{1k})-q(z_{ik}-
z_{i1}z_{11}^{-1}z_{1k})(z_{jk}-z_{j1}z_{11}^{-1}z_{1k})=  
$ 
 
$\underbrace{z_{jk}z_{ik}}_{\frame{1}} - \underbrace{z_{jk}z_{i1}
z_{11}^{-1}z_{1k}}_{\frame{2}} - \underbrace{z_{j1}z_{11}^{-1}z_{1k}
z_{ik}}_{\frame{3}} + 
\underbrace{z_{j1}z_{11}^{-1}z_{1k}z_{i1}z_{11}^{-1}z_{1k}}_{\frame{4}} +$ 
 
$ - \underbrace{qz_{ik}z_{jk}}_{\frame{1'}} + q \underbrace{z_{i1}
z_{11}^{-1}z_{1k}z_{jk}}_{\frame{2'}} + q \underbrace{z_{ik}z_{j1}
z_{11}^{-1}z_{1k}}_{\frame{3'}} - q \underbrace{z_{i1}z_{11}^{-1}
z_{1k}z_{j1}z_{11}^{-1}z_{1k}}_{\frame{4'}}$ 
 
Now $\frame{1}=\frame{1'}$ and $\frame{4}=\frame{4'}$ and: 
 
$\frame{2'}=z_{11}^{-1}z_{1k}z_{i1}z_{jk}$ 
 
$\frame{3}=q^{-1}z_{11}^{-1}z_{1k}z_{ik}z_{j1}$ 
 
$\frame{2}%
=q^{-1}(z_{jk}z_{11}^{-1})z_{i1}z_{1k}=q^{-1}z_{11}^{-1}z_{1k}
(z_{jk}z_{i1})-(q^{-1}-q)q^{-2}z_{11}^{-2}z_{1k}^2z_{i1}z_{j1}=  
$ 
 
$=z_{11}^{-1}z_{1k}z_{i1}z_{jk}\ 
-(q^{-1}-q)z_{11}^{-1}z_{1k}z_{ik}z_{j1}-(q^{-1}-q)q^{-2}
z_{11}^{-2}z_{1k}^2z_{i1}z_{j1}  
$ 
 
$\frame{3'}%
=(z_{ik}z_{11}^{-1})z_{j1}z_{1k}=qz_{11}^{-1}z_{1k}z_{ik}z_{j1}+
(q^{-1}-q)q^{-2}z_{11}^{-2}z_{1k}^2z_{i1}z_{j1}  
$  
 
\smallskip\  
 
The cancellations are obvious now . \qed
 
\begin{lemma} 
\label{lem-ad-da}
$AD-DA=(q^{-1}-q)BC$ 
\end{lemma} 
 
\pf\smallskip\  
 
$AD-DA-(q^{-1}-q)BC=$ 
 
$ 
(z_{ik}-z_{i1}z_{11}^{-1}z_{1k})(z_{jl}-z_{j1}z_{11}^{-1}z_{1l})
 - (z_{jl}-z_{j1}z_{11}^{-1}z_{1l})(z_{ik}-z_{i1}z_{11}^{-1}z_{1k})-  
$ 
 
$ 
-(q^{-1}-q)(z_{il}-z_{i1}z_{11}^{-1}z_{1l})(z_{jk}-z_{j1}z_{11}^{-1}z_{1k})=  
$ 
 
$\underbrace{z_{ik}z_{jl}}_{\frame{1}} -
\underbrace{z_{i1}z_{11}^{-1}z_{1k}z_{jl}}_{\frame{2}} -
\underbrace{z_{ik}z_{j1}z_{11}^{-1}z_{1l}}_{\frame{3}} +
\underbrace{z_{i1}z_{11}^{-1}z_{1k}z_{j1}z_{11}^{-1}z_{1l}}_{\frame{4}} +$ 
 
$-\underbrace{z_{jl}z_{ik}}_{\frame{1'}} 
+ \underbrace{z_{jl}z_{i1}z_{11}^{-1}z_{1k}}_{\frame{2'}} + 
\underbrace{z_{j1}z_{11}^{-1}z_{1l}z_{ik}}_{\frame{3'}} - 
\underbrace{z_{j1}z_{11}^{-1}z_{1l}z_{i1}z_{11}^{-1}z_{1k}}_{\frame{4'}}$ 
 
$-(q^{-1}-q)(\underbrace{z_{il}z_{jk}}_{\frame{1''}}
-\underbrace{z_{il}z_{j1}z_{11}^{-1}z_{1k}}_{\frame{2''}} -
\underbrace{z_{i1}z_{11}^{-1}z_{1l}z_{jk}}_{\frame{3''}}+
\underbrace{z_{i1}z_{11}^{-1}z_{1l}z_{j1}z_{11}^{-1}z_{1k}}_{\frame{4''}})$ 
 
Now $\frame{1}-\frame{1'}=\frame{1''}$ and $\frame{4}=q\frame{4'}$ , 
$\frame{%
4'}=q^{-1}\frame{4''}$ and we are left with: 
 
$AD-DA-(q^{-1}-q)BC=(\frame{2'}-(q^{-1}-q)\frame{2''})+(\frame{2}-(q^{-1}-q)  
\frame{3''})+$ $\frame{3'}-\frame{3}=$ 
 
$%
=(z_{jl}z_{i1}-(q^{-1}-q)z_{il}z_{j1})z_{11}^{-1}z_{1k}-z_{i1}
z_{11}^{-1}(z_{1k}z_{jl}-(q^{-1}-q)z_{1l}z_{jk})+  
\frame{3'}-\frame{3}=$ 
 
$=\underbrace{z_{i1}z_{jl}z_{11}^{-1}z_{1k}}_{\frame{5}} -%
\underbrace{ z_{i1}z_{11}^{-1}z_{jl}z_{1k}}_{\frame{6}} 
+ \frame{3'} - \frame{3}$ 
 
$\frame{5}%
=z_{i1}(z_{jl}z_{11}^{-1})z_{1k}=z_{i1}z_{11}^{-1}z_{jl}z_{1k}+
(q^{-1}-q)q^{-2}z_{i1}z_{11}^{-2}z_{1l}z_{j1}z_{1k}=  
$ 
 
$%
=q^{-1}z_{11}^{-1}z_{i1}(z_{1k}z_{jl})+(q^{-1}-q)q^{-3}z_{11}^{-2}
z_{1k}z_{1l}z_{j1}z_{i1}=  
$ 
 
$%
=q^{-1}z_{11}^{-1}z_{i1}z_{jl}z_{1k}-(q^{-1}-q)q^{-1}z_{11}^{-1}
z_{1l}z_{i1}z_{jk}+(q^{-1}-q)q^{-3}z_{11}^{-2}z_{1k}z_{1l}z_{j1}z_{i1}  
$ 
 
$\frame{6}%
=q^{-1}z_{11}^{-1}z_{i1}(z_{jl}z_{1k})=q^{-1}z_{11}^{-1}z_{i1}z_{1k}
z_{jl}-(q^{-1}-q)q^{-1}z_{11}^{-1}z_{1l}z_{i1}z_{jk}  
$ 
 
$\frame{3}%
=q^{-1}(z_{ik}z_{11}^{-1})z_{j1}z_{1l}=q^{-1}z_{11}^{-1}z_{ik}z_{1l}
z_{j1}-(q^{-1}-q)q^{-2}z_{11}^{-2}z_{1k}z_{i1}z_{1l}z_{j1}=  
$ 
 
$%
=q^{-1}z_{11}^{-1}z_{1l}z_{ik}z_{j1}-(q^{-1}-q)q^{-3}z_{11}^{-2}z_{1k}
z_{1l}z_{i1}z_{j1}  
$ 
 
$\frame{3'}=q^{-1}z_{11}^{-1}z_{1l}z_{ik}z_{j1}$ 

 \non This also ends the proof of the Lemma and the proof of the 
Theorem. \qed
\smallskip\  
 
\brk A warning is in order, it is not true that by any elementary 
transormation a quantum matrix would change to a new quantum matrix. 
We  proved that this is only in the case when one  performs the usual 
elementary transformations such that all elements in the first column 
except just one are set to zero, this is the matrix we denoed by $Z'$. 
Of course, now this can be iterated several times in order to obtain 
an upper triangular matrix.
\erk
\brk

{\bf The case $n\;  =\; 3$} 
\label{ex-n=3}
 \non The quantum determinant is: 
 $$ det_{q}(Z_3)=z_{11}z_{22}z_{33}-q^{-1}z_{12}z_{21}z_{33}-q^{-1}
z_{11}z_{23}z_{32}+q^{-2}z_{12}z_{23}z_{31}+q^{-2}z_{13}z_{21}
z_{23}-q^{-3}z_{13}z_{22}z_{31}.$$ 
\non One can prove through  laborious calculations {\em (see [PH])} that 
the coset of the $det_q(Z_3)$ mod commutators equals 
the Dieudonn\'e determinant of the $3$ by $3$ quantum matrix $Z_3$. 
 
\erk

\brk

At this point, induction would work if the quantum determinant were  
multiplicative. This is not the case. For instance for the quantum matrix 
$Z$ corresponding to the parameter $q$,  the square matrix $Z^2$ is a 
quantum matrix for the parameter $q^2$, this is a quite astonishing fact 
that calls for explanation!  Due to the fact that $det_q(Z)$ is a group-like
 element, one has  a weak multiplicative property 
namely if  $A$ and $B$ are quantum matrices such that their elements 
commute pairwise, then the product $AB$ is again a quantum matrix and 
$ det_q(AB) = det_q(A)det_q(B).$ 
 
But in the case when we perform elementary row operations on $Z$, say $A$ is 
a product of transvections and $B=Z$, the components of these matrices do 
not commute, so the above-mentioned result does not apply. 
\erk 

Meanwhile we can prove directly the following result which is all we need: 
 
\begin{theorem} [Row reducing the quantum determinant]
\label{thm-q-ind-step} 
The quantum determinants of the matrix $Z'$ obtained by row-reducing a 
quantum matrics $Z$ as in formulas (\ref{z'-1-st-row}) and (\ref{z'-i-th-row}) 
is equal to the quantum determiant of the original quantum matrix:
 \beqn \left| Z\right| _q=\left| Z^{\prime }\right| _q \eeqn 
\end{theorem} 
 
\pf  
Let $E^{\prime }=k<z_{ij},z_{11}^{-1}>,R^{\prime }=E^{\prime }\otimes _kB,$ 
and $\xi _i^{\prime }=\sum_jz_{ij}\xi _j,$ and $\xi _i^{\prime \prime 
}=\sum_jz_{ij}^{\prime }\xi _j$ . 
 
In $R^{\prime }$ we have: 
$$\xi _1^{\prime }\xi _2^{\prime }\cdots \xi _n^{\prime }=
det_q(Z)\ \xi _1\xi _2\cdots \xi _n\quad \mbox{\ and \ } \quad \xi _1^{\prime \prime }\xi _2^{\prime 
\prime }\cdots \xi _n^{\prime \prime }=det_q(Z^{\prime })\ \xi _1\xi 
_2\cdots \xi _n .$$ 
If we denote$\quad t_i:=-z_{i1}z_{11}^{-1}$ then $\xi _i^{\prime \prime 
}=\xi _i^{\prime }+t_i\cdot \xi _1^{\prime }$ for any $i>1$ we may write: 
 $$\qquad \xi _1^{\prime \prime }\xi _2^{\prime \prime }\cdots \xi _n^{\prime 
\prime }=\xi _1^{\prime }(\xi _2^{\prime }+t_2\xi _1^{\prime })\cdots (\xi 
_n^{\prime }+t_n\xi _1^{\prime }).$$ 
 By induction we prove now: 
\beqn \xi _1^{\prime }(\xi _2^{\prime }+t_2\xi _1^{\prime })\cdots (\xi 
_n^{\prime }+t_n\xi _1^{\prime })=\xi _1^{\prime }\xi _2^{\prime }\cdots \xi 
_{i-1}^{\prime }(\xi _i^{\prime }+t_i\xi _1^{\prime })\cdots (\xi _n^{\prime 
}+t_n\xi _1^{\prime })\eeqn 
 Indeed: 
\beqa \xi _1^{\prime }\xi _2^{\prime }\cdots \xi _{i-1}^{\prime }
(\xi _i^{\prime 
}+t_i\xi _1^{\prime })\cdots (\xi _n^{\prime }+t_n\xi _1^{\prime })&=&\\ 
\xi _1^{\prime }\xi _2^{\prime }\cdots \xi _{i-1}^{\prime }\xi 
_i^{\prime }(\xi _{i+1}^{\prime }+t_{i+1}\xi _1^{\prime })\cdots (\xi 
_n^{\prime }+t_n\xi _1^{\prime })&+&\xi _1^{\prime }\xi _2^{\prime }\cdots \xi 
_{i-1}^{\prime }\cdot t_i\xi _1^{\prime }(\xi _{i+1}^{\prime }+t_{i+1}\xi 
_1^{\prime })\cdots (\xi _n^{\prime }+t_n\xi _1^{\prime })= \\ 
 \xi _1^{\prime }\xi _2^{\prime }\cdots \xi _{i-1}^{\prime }\xi 
_i^{\prime }(\xi _{i+1}^{\prime }+t_{i+1}\xi _1^{\prime })\cdots (\xi 
_n^{\prime }+t_n\xi _1^{\prime })&+&\xi _1^{\prime }\xi _2^{\prime }\cdots \xi 
_{i-1}^{\prime }\cdot \xi _1^{\prime }qt_i(\xi _{i+1}^{\prime }+t_{i+1}\xi 
_1^{\prime })\cdots (\xi _n^{\prime }+t_n\xi _1^{\prime })
\eeqa 
\noindent because a direct calculation shows that $t_i\xi _1^{\prime }=q\xi 
_1^{\prime }t_i.$ But using the fact that $\xi _i^{\prime }$ are points of 
the quantum Grassmannian plane, by Proposition 2 i.e., $\xi _i^{\prime 2}=0$ 
for all $i$, and $\xi _i^{\prime }\xi _j^{\prime }=q^{-1}\xi _j^{\prime }\xi 
_i^{\prime }$ for all $i<j$, it follows that $\xi $ $_1^{\prime }\xi 
_2^{\prime }\cdots \xi _{i-1}^{\prime }\cdot \xi _1^{\prime }=0$ so the 
second term vanishes. We shall later need a multiparameter version of this 
result, you may see it holds as well. Eventually for $i=n$ we get: 
$$\xi _1^{\prime \prime }\xi _2^{\prime \prime }\cdots \xi 
_n^{\prime \prime }=\xi _1^{\prime }(\xi _2^{\prime }+t_2\xi _1^{\prime 
})\cdots (\xi _n^{\prime }+t_n\xi _1^{\prime })=\xi _1^{\prime }\xi 
_2^{\prime }\cdots \xi _n^{\prime }$$
\noindent therefore
$$\xi _1^{\prime }\xi _2^{\prime }\cdots \xi _n^{\prime }=det_q(Z)\ \xi _1\xi 
_2\cdots \xi _n=\xi _1^{\prime \prime }\xi _2^{\prime \prime }\cdots \xi 
_n^{\prime \prime }=det_q(Z^{\prime })\ \xi _1\xi _2\cdots \xi _n$$ 

\noindent hence: $det_q(Z)$ =$det_q(Z^{\prime }).$ \qed
 
\smallskip\  
 
\noindent We can now prove 
 
\becor [Quantum determinant and Dieudonne\'e determinant]
\label{cor-qdet=ddet}
Let $Z$ be a quantum matrix i.e., its elements satisfy the relations 
$\ {\cal R%
}el_q.$ Then 
 
\qquad \qquad  
$$ 
\pi (det_q(Z))=det_D(Z).  
$$ 
 \eecor
\pf  
\noindent Using the result above  above 
$\left| Z\right| _q=\left| Z^{\prime 
}\right| _q$. Now by column expansion $\left| Z^{\prime }\right| 
_q=z_{11}\left| Z^{\prime \prime }\right| _q.$ 
 
\noindent Because $Z^{\prime }$ hence also $Z^{\prime \prime }$ are quantum 
matrices and the dimension of $Z^{\prime \prime }$ is $n-1,$ by induction 
$ \pi (\left| Z^{\prime \prime }\right| _q)=\det _D(Z^{\prime \prime }).$ 
On the other hand by Theorem \ref{thm-well-def-ddet} 
$det_D(Z)=\pi (z_{11})det_D(Z^{\prime \prime }).$  \qed

\section{\bf Row-reducing in the  multiparametric case} 
 
 
In this section we look at the quantum determinant for the multiparameter 
quantum linear group over a field k. We shall use the notation and results 
from [AST]. 
 
If $\ p=(p_{ij})$ is such that $p_{ij}p_{ji}=1$ and $p_{ii}=1,\ $we call 
such a matrix $p$ a multiplicatively antisymmetric matrix. If $\lambda 
,p_{ij}$ $\in k^{*}$ define the multiparameter quantum linear group to be a 
k-algebra with generators $u_{ij},\ M_{p,\ \lambda 
}(n)=k<u_{11,}u_{12,}\cdots ,u_{nn}>$ subject to the following  
{\it twisted quantum relations}: 
  
\beqa TWRel_{p,\lambda}\left\{ \begin{array}{l}
  \mbox{\it Row relations:\ } u_{il}\ u_{ik}=\frac 1{p_{lk}\ }\ 
u_{ik}\ u_{il} \mbox {\it \  for  all \ } k<l  \label{mpmq-rel-row} \\
  \mbox{\it Column relations:\ } u_{jk}\ u_{ik}=\lambda \ p_{ji}\ 
u_{ik}\ u_{jk}
\mbox {\it  \ for all \  }  i<j \label{mpmq-rel-col} \\
   \mbox{\it Secondary diagonal relations:\ } u_{jk}\ u_{il}=
\frac{\lambda p_{ji}}{ p_{kl}\ }\ u_{il}\ u_{jk} \mbox {\it \ 
for all \ } i<j \mbox {\it ,} k<l \label{mpmq-rel-sdiag} \\
  \mbox{\it Main diagonal relations:\ } u_{jl}\ u_{ik}=
\frac{p_{ji}}{p_{kl}}\ 
u_{ik}\ u_{jl}+(\lambda -1) p_{ji}\ u_{il}\ u_{jk}\  
\mbox {\it if \ }  i<j \mbox {\it ,} k<l . \label{mpmq-rel-mdiag} 
\end{array}
\right.
\eeqa
 
\noindent We denote by $U$ the matrix $(u_{ij})$, with $u_{ij}$ satisfying 
the above twisted quantum relations $TWRel_{p,\lambda }$ , and let $%
q_{ji}:=\lambda p_{ji}$ for all  $i<j$ , $q_{ii}:=1$ and $q_{ji}:=\lambda 
^{-1}p_{ij}$ for $i>j$, and $I=\{1,\cdots ,n\}.$ 
 
\noindent The comultiplication on $M_{p,\lambda }(n)$ is dual to matrix 
multiplication:$\ \Delta (u_{ij})=\sum\limits_ku_{ik}\otimes u_{kj},$ and 
the counit is defined by $\epsilon (z_{ij})=\delta _{ij,}$ both extended 
multiplicatively to all $M_{p,\lambda }(n).\ $Then $\ M_{p,\lambda }(n)$ \ 
has a bialgebra structure.$\ M_{p,\lambda }(n)$ can be obtained by using the 
Manin quadratic algebra construction for the following two quadratic 
algebras: 
 
 \begin{itemize} 
\item  $A=k<x_1,\cdots ,x_n>\ $subject to relations : $x_jx_i=q_{ji}x_ix_j$
 for all  
$\ i,j\in I$ 
 
\item $B=k<y_1,\cdots ,y_n>\ $subject to relations : $y_jy_i=p_{ij}y_iy_j\ $
for 
all $\ i,j\in I$ 
 
\end {itemize} 
 
We proceed as in [AST] to define the determinant, namely we consider the 
corresponding Grassmannian algebras $A^{!}$ and $B^{!}$ 
 
 \begin{itemize}  
\item $A^{!}=k\negthinspace <\negthinspace \xi _1,\cdots ,
\xi _n\negthinspace >\ $%
subject to relations: $\xi $$_j\xi _i=-q_{ij}\xi _i\xi _j$ and $\xi 
_i^2=0,\forall $ $i,j\in I$ 
 
\item $B^{!}=k\negthinspace \negthinspace <\negthinspace \eta _1,\cdots ,
\eta _n%
\negthinspace >\ \ $subject to relations: $\eta $$_j\eta _i=-p_{ji}\eta 
_i\eta _j\ $and $\eta $$_i^2=0,\forall $ $i,j\in I$ 
  \end{itemize} 
\smallskip\  
 
\noindent For any subset $J$ of $I$ write $\xi _J=\prod_{j\in J}\xi _j$ 
and $%
\eta $$_J=\prod_{j\in J}\eta _j$ the product is taken in increasing order. 
 
\noindent Let $H=H(p,\lambda )=M_{p,\lambda }(n),$ and consider the coactions: 
 
$\delta :A^{!}\longrightarrow A^{!} \otimes H $ and $\delta 
:B^{!}\longrightarrow H\otimes B^{!}$ they are both algebra homomorphisms. 
 
\noindent If $\delta ($ $\xi _j)=$ $\sum_ju_{jk}\otimes \xi _j$ and 
$\delta ($ $\eta $$%
_j)=$ $\sum_ju_{jk}\otimes \eta _j$ \ then if follows that 
 
$\delta ($ $\xi _K)=\sum_J \xi _J \otimes U_{JK}$ and $\delta ($ $\eta 
_K)=\sum_JU_{KJ}\otimes \eta _J$ 
 
\noindent where, as in [AST], Lemma 1 : $[J,K]$ is the set of bijective 
mappings $\theta $$:J\widetilde{\longrightarrow }K$ for any subsets $J,K$ of  
$I$ and by definition
 $$U_{J,K}:=\sum_{\theta \in [J,K]}\sigma (p,\theta )\prod_{j\in J}
u_{j,\theta 
j}=\sum_{\theta \in [J,K]}\sigma (q^{-1},\theta )\prod_{k\in K}u_{\theta 
k,k} $$  
with the products taken in increasing order and $\sigma (p,\theta 
)=\prod_{ j<j^{\prime }, \: \theta j>\theta j^{\prime }} (-p_{\theta 
j,\theta j^{\prime }}).$ 
 
In particular for $J=K=I$\ we get$\ det_{p,\lambda }(U)=U_{I,I}.$ 
 
\begin{definition} 
\label{def-multipar-qdet}
The multiparameter quantum determinant is: 
 
\beqn
det_{p,\lambda }(U)=U_{I,I}=\sum_{\theta \in S_n}\sigma (p,\theta )\ 
u_{1,\theta 1}\cdots u_{n,\theta n} \mbox{ \it \ \ were   \ } 
\sigma (p,\theta 
)=\prod_{ j<j^{\prime }, \: \theta j>\theta j^{\prime }} 
(-p_{\theta 0j,\theta j^{\prime }})
\eeqn
\end{definition} 
 
\  
 
As we did before for the quantum matrices, if $H=k<u_{ij}>,R=H\otimes _kB,$ 
then a reformulation of this definition is that in $R$ we have: 
 
\begin{proposition}
\label{prop-multipa-Kobyzev} 
If $\eta _i^{\prime }=\sum_ju_{ij}\eta _j\ $ then$\quad \eta _1^{\prime }\eta 
_2^{\prime }\cdots \eta _n^{\prime }=det_{p,\lambda }(U)\ \eta _1\eta _2\cdots 
\eta _n.$ 
\end{proposition} 
 
In a similar way we have formulas for rows and columns expansion. Denote by $%
U_{ji}$ the determinant $det_{p,\lambda }$ of the $(n-1)$ by $(n-1)$ matrix 
obtained by removing the $j$-th row and the $i$-th column of $U.$ \smallskip%
\  
 
Let us denote $\beta _j:=\prod_{m=j+1}^n(-q_{jm})$ and
 $\gamma _j:=\prod_{m=1}^{j-1}(-p_{jm}).$ 
 
The following is proved in [AST] Thm 3, (21). 
 
\begin{proposition} 
\label{prop-multipar-cofactor}
If U is a quantum multiparameter matrix then 
 
1) The element $det_{p,\lambda }(U)\ $ is normalizing (but it is not central) 
 
2) The following row and column expansion formulas hold: 
\end{proposition} 
 
\beqn 
det_{p,\lambda }(U)\ \ =\sum_j\frac{\beta _j}{\beta _k}U_{jk}\ u_{jk}=\sum_k  
\frac{\gamma _k}{\gamma _j}\ u_{jk}U_{jk}.  
\eeqn 
 
Another important fact established in [AST] that we need to use is the fact 
that $\ M_{p,\lambda }(n)$ can be obtained by twisting the multiplication in  
$M_q(n)$ by a cocycle associated to $p$. 
 
First note that for $\lambda =q^2$ and $p_{ij}^{\prime }=q$ , for $i<j,\ 
p_{ij}^{\prime }=q^{-1}$ for $i>j$ and $p_{ii}^{\prime }=1$ the 
corresponding multiparametric quantum matrix ring $M_{p^{\prime },q^2}(n)$ 
is just$\ M_q(n).$ 
 
Now let $G\ $be an abelian group isomorphic to the product of n copies of $Z$%
. Using the multiplicative notation and denoting by $t_i\ $ the generator of 
the $i$-th copy we may write 
\beqn \qquad G\approx \ <t_1>\times <t_2>\times \cdots \times <t_n>.\eeqn 
We give $M_{p,\lambda }(n)$ a grading by letting $u_{ij}$ have left degree 
equal to $t_i$ and right degree equal to $t_j$ and extend this 
multiplicatively. Then if we remark that $M_{p,\lambda }(n)$ is an iterated 
Ore extension, and if we let $D$ be its total field of fractions, then $%
u_{ij}^{-1}$, an element in $D$, has left degree $t_i^{-1}$and right degree
 $t_j^{-1}.$ 
 
For an arbitrary antisymmetric matrix $p=(p_{ij})$ such that $p_{ij} p_{ji}=1$ 
and $p_{ii}=1$ define the $2$-cocycle on $G$ by: 
  $$c_p(\prod t_i^{m_i},\prod t_j^{n_j}):=\prod_{i<j}\ p_{ij}^{m_in_j},$$ 
\non in fact this is the unique bimultiplicative function such that: 
$$ c_p(t_i,t_j) = p_{ij} \mbox{\ if\ } i<j 
\mbox{\  and \ } c_p(t_i,t_j) = 1 \mbox{\  if\ } i\geq j.$$ 
 Conversely for any $2$-cocycle $c$ on $G$ we associate the following 
matrix $r(c),$ which has the property that $p$ is antisymmetric: 
$$ r_{ij}(t_i,t_j)=\frac{c(t_i,t_j)}{c(t_j,t_i)}.$$ 
 
 
Proposition 1 in [AST ] shows that the correspondences defined above 
define a bijection from $H^2(G,k^{*})$ to the set of multiplicatively 
antisymmetric matrices $p$, i.e., with the property $p_{ij}p_{ji}=1$ and $%
p_{ii}=1$ 
 
If $A$ is a $G\times G$-graded algebra, having both a left and a right 
grading we define a new multiplication ''$\circ $ '' on A called the 
multiplication twisted by $c_p{}^{-1}$ on the left and by $c_p$ on the right 
by the formula:%
\beqn  
a\circ b=c_p^{-1}(t_i,t_j)\ c_p(t_k,t_l)\ a\cdot b  
\eeqn \label{def-twist-multipl}
 \noindent where $a$ has left degree $t_i$ and right $\deg $ree $t_k$, $b$ 
has left degree $t_j$ and right $\deg $ree $t_l$. 
 
We write $_{c^{-1}}A_c$ for the algebra $A$ with the twisted multiplication $%
\circ .$ 
 
We need also the following result established in [AST] Thm 4: 
 
\begin{proposition} 
\label{prop-twist-multipar}
If we twist $H(\ p,\lambda ) = M_{p,\ \lambda }(n)$ simultaneously by $%
c_p{}^{-1}$ on the left and by $c_p$ on the right we obtain 
 
$$ 
{}_{c^{-1}} M_{p,\ \lambda }(n)_c = M_{r(c)p,\lambda }(n) 
\mbox{ \it \ \ or\ \ \  }{}_{c^{-1}} H(\ p, \lambda )_c = H( r(c)p, 
\lambda) .  
$$ 
\end{proposition} 
  
\noindent In particular $M_{p,\ \lambda }(n)$ can be obtained by twisting $%
M_q(n)$ = $H(p^{\prime },\lambda )$ (for $\lambda =q^2$ and $p_{ij}^{\prime 
}=q$ , for $i<j,\ p_{ij}^{\prime }=q^{-1}$ for $i>j$ and $p_{ii}^{\prime }=1$%
) by the cocycle defined above $c_p.$ In this particular case we shall not 
use $\circ $ any more for the multiplication but simply juxtaposition. 
 
\smallskip\  
 
\noindent We are now ready now to look at the determinant of the matrix $U$. 
First we look at 
 
\begin{example}[ The multiparametric quantum determinant for $n \: = \: 2$]
\label{ex-multipar-qdet-n=2} 
 \end{example} 
 
The multiparametric quantum determinant is: 
  
$$\left|  
\begin{array}{cc} 
u_{11} & u_{12} \\  
u_{21} & u_{22} 
\end{array} 
\right| _{p,\lambda }=u_{11}u_{22}-p_{21}u_{12}u_{21}$$ 
 
The Dieudonn\'e determinant is 
  
$det_D(U)=u_{11}(u_{22}-u_{21}u_{11}^{-1}u_{12}).$ 
  
Now use the defining relations: $u_{21}u_{11}^{-1}=\frac 
1{\lambda p_{21}}u_{11}^{-1}u_{21}$ and $\ u_{21}u_{12}=
\frac{\lambda p_{21}}{p_{12}}%
\ u_{12}u_{21}.$

Use also $\ \frac 1{p_{12}}=p_{21}$ and this will establish $%
\pi (det_{p,\lambda }(U))=det_D(U)$ for $n=2.$ 
 
\smallskip\  
 
We may now prove that if $U$ is a multiparametric quantum matrix satisfying 
the twisted relations\ $TWRel_{p,\lambda }$ then for any $n$ we have%
$\ \pi (det_{p,\lambda }(U))=det_D(U).$

 Let us define the matrix $U^{\prime }=(u_{ij}^{\prime })$ by  
 
\beqn 
u_{1j}^{\prime }=u_{1j} \mbox{\it \ for any \ } j,u_{i1}^{\prime }=0,
 \mbox{ \it \ for\ }  i=2,\ldots ,n, \label{u'1}\\ 
u_{ij}^{\prime }=u_{ij}-u_{i1}u_{11}^{-1}u_{1j} \mbox{\it \ for \ } 
i=2,\ldots ,n \mbox{\it \ and \ } j=2,\ldots ,n.\label{u'i}
\eeqn 
$U^{\prime }$ is obtained from $U\ $ by clearing up all the positions of the 
first column except the first one by use of elementary row transformations. 
 
Let $U^{\prime \prime }$ be the $(n-1)\times (n-1)$ matrix obtained by 
deleting the first row and the first column of $U^{\prime }.$ Just like 
before we can establish the following
 
\begin{theorem} [Row-reducing the multiparametric quantum matrix]
\label{thm-row-red-multipar}
The matrix $U^{\prime }$ is a multiparametric quantum matrix, i.e., its 
elements satisfy the relations $TWRel_{p,\lambda }$ defining $%
M_{p,\lambda }(n)$ 
\end{theorem} 
 \pf 
The proof relies on the following fact: each relation in $TWRel_{p,\lambda }$
 is obtained by twisting the corresponding relation in 
$ Rel_q$ (the relations of $M_q(n)).$ We use the identification
 $z_{ij}\longleftrightarrow u_{ij}$ and keep in mind that the product of 
the $u_{ij}$ is obtained by twisting by $c_p{}^{-1}$ on the left and by $c_p$ 
on the right the multiplication of the corresponding $z_{ij}$, so we get a 
factor like $c_p(t_i,t_j)\ c_p^{-1}(t_k,t_l)$ when twisting the product 
$z_{ik}\cdot z_{jl}.$ 
 
Now a direct calculation shows that $u_{i1}u_{11}^{-1}u_{1j}$ comes from the 
corresponding product $z_{i1}z_{11}^{-1}z_{1j}$ changed by the factor 
$c_p(t_{i,}\ t_1^{-1})$ $c_p^{-1}(t_{1,}\ t_1^{-1})\ 
c_p(t_it_1^{-1},t_1^{-1})\ c_p^{-1}(e_{,}\ t_j).$
 This factor is equal to $1$ 
because $c_p$ is bimultiplicative by its definition. It means the product is 
not changed by the twist. It has the same left and same right degree as
 $u_{ij}.\;$ 
 
Therefore $u_{ij}^{\prime }$ has left degree $t_i$ and right degree $t_j,$ 
in fact the same left and same right degree as $u_{ij}.\;$Then by twisting ( 
by c$_p^{-1}$ on the left and by c$_p$ on the right) the products in the 
relations $Rel_q$ satisfied by $z_{ij}^{\prime }$ we get that the 
relations $TWRel_{p,\lambda }$ are satisfied by $u_{ij}^{\prime },$ 
because the $\deg $rees are the same for like for $z_{ij}$ and $%
z_{ij}^{\prime }.$ We can see now that $u_{ij}^{\prime }$ satisfy exactly 
the twisted relations $TWRel_{p,\lambda }.$ This establishes the 
theorem. \qed
 
\becor [Multiparametric quantum determinant and Dieudonne\'e determinat]
\label{cor-multipar-qdet=ddet}
Let $U$ \ be a multiparametric quantum matrix satisfying the relations\  
$TWRel_{p,\lambda }.$ Then for any $n$ we have  
$$ 
\pi (det_{p,\lambda }(U))=det_D(U).  
$$ 
\eecor 
  \pf  
\noindent We use induction exactly as in the proof of the quantum determinant.

\noindent The proof for $\left| Z\right| _q=\left| Z^{\prime }\right| _q$ 
works exactly in the same way, establish that $\left| U\right| _{p,\lambda 
}=\left| U^{\prime }\right| _{p,\lambda }$ and by column expansion $\left| 
U^{\prime }\right| _q=u_{11}\left| U^{\prime \prime }\right| _q.$ 
 
Because $U^{\prime }$ hence also $U^{\prime \prime }$ are multiparametric 
quantum matrices and the dimension of $U^{\prime \prime }$ is $n-1,$ by 
induction $\pi ($$\left| U^{\prime \prime }\right| _q)=\det _D(U^{\prime 
\prime }).$ On the other hand by Theorem \ref{thm-well-def-ddet}, 
$det_D(U)=\pi 
(u_{11})det_D(U^{\prime \prime })$ which proves our result. \qed
 
\brk 
We proved  the Theorems \ref{thm-row-red-qmatr}  and
 \ref{thm-row-red-multipar} and their corollaries  for the generic 
matrices $Z$ and $U.$ In this 
case we do not have to worry about the existence of $\ z_{11}^{-1}$ or $%
u_{11}^{-1}.$ For an arbitrary quantum matrix or quantum multiparameter 
matrix we may now use specialization. 
\erk

\end{document}